\documentclass[a4paper,11pt]{scrartcl}

\usepackage[utf8]{inputenc}
\usepackage{amsmath, amsthm, amssymb}
\usepackage{xcolor,soul}
\usepackage{xspace}
\usepackage{mathtools}
\usepackage{rotating}

\usepackage[none]{optional}
\newcommand{\plan}[1]{}
\newcommand{\ba}[1]{}
\newcommand{\pn}[1]{}
\opt{draft}{
\renewcommand{\plan}[1]{\textcolor{blue}{#1}}
\renewcommand{\ba}[1]{\textcolor{orange}{BA: #1}}
\renewcommand{\pn}[1]{\textcolor{purple}{PN: #1}}
}

\usepackage[color,all,2cell]{xy}\SelectTips{cm}{12}\SilentMatrices\UseAllTwocells  %
\usepackage{hyperref}
\hypersetup{
    hidelinks
}

\theoremstyle{plain}
\newtheorem{thm}{Theorem}

\theoremstyle{remark}

\newtheorem{exercise}[thm]{Exercise}

\theoremstyle{definition}

\newcommand{\Type}{\constfont{Type}}
\newcommand{\Set}{\constfont{Set}}
\newcommand{\Prop}{\constfont{Prop}}

\def\lv{\mathopen{{[\kern-0.14em[}}}    %
\def\rv{\mathclose{{]\kern-0.14em]}}}   %

\newcommand{\refl}{\constfont{refl}}

\newcommand{\constfont}[1]{\ensuremath{\mathsf{#1}}}
\newcommand{\C}{\ensuremath{\mathcal{C}}}

\newcommand{\transport}{\ensuremath{\mathsf{transport}}}
\newcommand{\isEquiv}{\ensuremath{\mathsf{isEquiv}}}
\newcommand{\transfib}[3]{\ensuremath{\transport^{#1}(#2,#3)\xspace}}
\newcommand{\id}[3][]{\ensuremath{#2 =_{#1} #3}\xspace}

\newcommand{\idtoequiv}{\constfont{idtoequiv}}
\newcommand{\idtoiso}{\constfont{idtoiso}}

\title{Univalent foundations \\ and the equivalence principle}
\author{Benedikt Ahrens \and Paige Randall North}
\date{}

\begin{document}

\maketitle

\begin{abstract}
In this paper, we explore the `equivalence principle' (EP): roughly,
statements about mathematical objects should be invariant under an appropriate notion of
equivalence for the kinds of objects under consideration.
In set theoretic foundations, EP may not always hold:
for instance, the statement `$1 \in \mathbb{N}$' is not invariant under isomorphism of sets.
In univalent foundations, on the other hand, EP
has been proven for many mathematical structures.
We first give an overview of earlier attempts
at designing foundations that satisfy EP.
We then describe how univalent foundations validates EP.
\end{abstract}

\section{The equivalence principle}

What should it mean for two objects $x$ and $y$ to be equal? One proposal by Leibniz \cite{leibniz}, known as the ``identity of indiscernibles'', 
states that if $x$ and $y$ have the same properties, then they must be equal:
 \begin{equation*}
       \forall \text{ properties } P, \left(P(x) \leftrightarrow P(y)\right) \enspace \to \enspace x = y  .
 \end{equation*}

For this proposal to be reasonable, then the converse, the ``indiscernibility of identicals,'' should hold incontrovertibly. That is, if $x$ and $y$ are equal, then they must have the same properties:

\begin{equation}
             x = y \enspace \to \enspace  \forall \text{ properties } P, \left(P(x) \leftrightarrow P(y)\right) .
       \label{eq:indiscernability_of_identicals}
\end{equation}

Indeed, one would be hard-pressed to find a mathematician who disagreed with this principle. However, in classical mathematics based on set theory, this principle is of limited usefulness:
too few objects are equal. A group theorist, for example, would have little interest in a principle which required them to suppose that two groups are equal.

Instead, mathematicians are often interested in weaker notions of sameness and those properties that are invariant under such notions. A group theorist, for example, would have more interest in an analogous principle that described the properties of any pair of isomorphic groups $G$ and $H$:
\begin{equation*}
             G \cong H \enspace \to \enspace  \forall \text{ group theoretic properties } P, \left(P(G) \leftrightarrow P(H)\right).  \label{eq:indiscernability_of_identical_groups}
             \end{equation*}
Similarly, category theorists would be more interested in a principle that described the properties of any pair of equivalent categories $A$ and $B$:
\begin{equation*}
             A \simeq B \enspace \to \enspace  \forall \text{ category theoretic properties } P, \left(P(A) \leftrightarrow P(B)\right). 
             \end{equation*}

To generalize: mathematicians working in some domain $\mathcal D$ often utilize a stronger variant of 
the principle given in line~\eqref{eq:indiscernability_of_identicals} above, called
the \textbf{equivalence principle}: for all objects $x$ and $y$ of domain $\mathcal{D}$:
 \begin{equation}    
       x \sim_{\mathcal{D}} y \enspace \to \enspace \forall~ \mathcal{D}\text{-properties } P, \left(P(x) \leftrightarrow P(y)\right)  \enspace ,
       \label{eq:equivalence_principle}
 \end{equation}
where \( \sim_{\mathcal{D}} \) denotes a suitable 
notion of sameness for the domain $\mathcal{D}$.

We might consider a still stronger variant of the equivalence principle. A group theorist, for example, might not only want \emph{properties of} groups to be invariant under isomorphism, but they might also want \emph{structures on} groups to be invariant under isomorphism. 
For example, if the equivalence principle \eqref{eq:equivalence_principle} holds in the domain of group theory and if two groups $G$ and $H$ are isomorphic, then the statements ``$G$ has a representation on $V$" and ``$H$ has a representation on $V$" are equivalent (for some fixed vector space $V$). However, it is actually the case that the isomorphism $G \cong H$ induces a bijection between the set of representations of $G$ on $V$ and the set of representations of $H$ on $V$, which we regard as structures on $G$ and $H$ respectively. Such a variant of the equivalence principle has become known as the \emph{Structure Identity Principle} (see \cite{oberwolfach},\cite[Section~9.8]{hottbook}, \cite{10.1093/philmat/nkt030}).

Our goal in this paper is to describe how one can find the right notion $ \sim_{\mathcal{D}} $ of sameness and the right class of `$\mathcal{D}$-properties and $\mathcal{D}$-structures' for some specific domains $\mathcal{D}$.

This right notion of sameness is not uniformly defined across different mathematical objects. However, we usually use the one already present in mathematical practice since we aim for the equivalence principle to capture mathematical practice.
As a rule of thumb, it is usually considered to be 
\begin{itemize}
 \item \emph{equality} when the objects naturally form a set---numbers,
 functions, etc.
 \item \emph{isomorphism} when the objects naturally form a category---sets, groups, etc.
 \item \emph{equivalence} when the objects naturally form a bicategory---e.g., categories.
\end{itemize}

The hard part will be in determining the right class of $\mathcal{D}$-properties and $\mathcal{D}$-structures for some specific domain $\mathcal{D}$.
In usual mathematical practice, we can state properties which break the equivalence principle; 
that is, we can state properties of mathematical objects that are not invariant under sameness. We will seek to exclude such properties from our class of $\mathcal{D}$-properties and $\mathcal{D}$-structures.

\begin{exercise}\label{ex:statement_variant_set}
  Denote by $2\mathbb{N}$ the set of even natural numbers.
  Find a property of sets that is not invariant under the isomorphism $\mathbb{N} \cong 2\mathbb{N}$ given by multiplying
  and dividing by $2$, respectively.
  \flushright \rotatebox[origin=c]{180}{\textbf{Answer}: One such statement is given in the abstract.}
\end{exercise}

\begin{exercise}\label{ex:equiv_singleton_cat_bool}
  Find a property of categories that is true for one, but not for the other
  of these two, equivalent, categories.
   
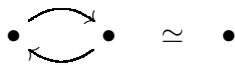
\[
     \xymatrix{ \bullet \rtwocell<4>{'} & \bullet} \quad \simeq \quad \bullet
   \]
   
   \flushright \rotatebox[origin=c]{180}{\textbf{Answer}: The statement ``The category $\C$ has exactly one object.'' is such a statement.}
\end{exercise}

Thus, to assert an equivalence principle for sets or categories, we need to exclude these properties from our collection of `set theoretic properties' and `category theoretic properties'.
M.\ Makkai \cite{MR1678360} says
\begin{quote}
    The basic character of the Principle of Isomorphism is that of a constraint on the
    language of Abstract Mathematics; a welcome one, since it provides for the separation of sense from nonsense.
\end{quote}
Put differently, establishing an equivalence principle means establishing a \emph{syntactic criterion} for properties and structures that are invariant under sameness.

\section{History} \label{sec:history}
Look again at Example~\ref{ex:equiv_singleton_cat_bool}. There, we violated the equivalence principle for categories by referring to equality of objects. This might lead one to conjecture (correctly) that categorical properties which obey the equivalence principle cannot mention equality of objects. 

However, the traditional definition of category mentions equality of objects. It usually includes the following axiom: for any two morphisms $f$ and $g$ such that the codomain of $f$ equals the domain of $g$, there is a morphism $gf$ such that domain of $gf$ equals the domain of $f$ and the codomain of $gf$ equals the codomain of $g$.

To avoid mentioning equality of objects, one can express the composability of morphisms of that category via different means,
specifically by having not one collection of morphisms but many \emph{hom-sets}: one for each pair of objects.
This idea, for instance explained in \cite[Section~I.8]{maclane} usually requires asking the hom-sets to be disjoint.
This last requirement is automatic if we work instead in a typed language, where types are automatically disjoint.

A category is then given by
  \begin{itemize}
  \item a type $O$ of objects,
  \item for each $x,y : O$, a type $A(x,y)$ of arrows from $x$ to $y$,
  \item for each $x,y,z : O$ and $f : A(x,y)$, $g : A(y,z)$, a composite arrow $g \circ f : A(x,z)$, and
  \item for each $x : O$, an identity arrow $\mathsf{id}_x : A(x,x)$
  \end{itemize}
such that
  \begin{itemize}
  \item for each $w,x,y,z : O$ and $f : A(w,x)$, $g : A(x,y)$, $h : A(y,z)$, there is an equality $h \circ (g \circ f) = (h \circ g) \circ f$ in $A(w,z)$,
  \item for each $x,y : O$ and $f : A(x,y)$, there is an equality $f\circ \mathsf{id}_x = f$ in $A(x,y)$, and
  \item for each $x,y : O$ and $f : A(x,y)$, there is an equality $ \mathsf{id}_y \circ f = f$ in $A(x,y)$.
  \end{itemize}
Note that when stating axioms, the only equality that is mentioned is the equality within a hom-set of the form $A(x,y)$,
that is, between arrows of the same type.

By adding quantifiers, ranging over one type at a time, to this typed language, we obtain a language for stating properties of, and constructions on, categories.
It turns out that the statements of that language are invariant under equivalence of categories:

\begin{thm}[\emph{Théorème de préservation par équivalence} \cite{MR539867}]
    \label{thm:blanc}
    A {property} of categories (expressed in 2-typed first order logic) is invariant under equivalence if and only if it can be expressed 
    in the typed language sketched above, and without referring to equality of objects.
\end{thm}
We do not give here the precise form of the typed language,  but refer instead to Blanc's article for details. 
Note that Freyd \cite{MR0412249} states a similar result to Blanc's above, in terms of ``diagrammatic properties''.

Makkai \cite{MR1678360} develops notions of \emph{signature} and \emph{theory}, to specify
mathematical structures.
A theory is a pair $(L,\Sigma)$ consisting of a signature $L$ (specifying the shape of the structure) and a set $\Sigma$ of ``axioms'' over $L$ (specifying the axioms of the structure).
A theory determines a notion of ``model''---which is an $L$-structure satisfying the properties specified by $\Sigma$---and
of ``equivalence'' of such models, called $L$-equivalence.

His \emph{Invariance Theorem} gives a result similar to Theorem~\ref{thm:blanc} for models of a theory:
given an interpretation $T = (L, \Sigma) \to S$ of such a 
theory in a first-order logic theory, an $S$-sentence $\phi$ is invariant under $L$-equivalence if and
only if it is expressible in First Order Logic with Dependent Sorts (FOLDS) over $L$.

In the following sections, we will see that similar results can be shown in univalent foundations.
Specifically, not only properties but also constructions will be ``invariant'' under equivalence,
and invariance of properties will be recovered as a special case via the propositions-as-some-types
correspondence.

\section{Univalent foundations and transport of structures along equivalences}

Starting in the 1970's, Per Martin-L\"of designed several versions of dependent type theory, which are now called Martin-L\"of Type Theories \cite{Martin-Lof-1972}.
These were intended to be foundations of mathematics that, unlike set theory, have an inherent notion of computation built in.
For decades, Martin-L\"of type theories have formed the basis of computer proof assistants such as Coq and Agda.

One of the most mysterious features of this kind of type theory is its \emph{equality} type $a =_X b$ of any two inhabitants $a$ and $b$ of a type $X$---see \cite[Section~4.3]{Altenkirch}.
Inhabitants of such an equality type behave, in many ways, like a proof of equality; 
in particular, they can be composed and inverted, corresponding to the transitivity and symmetry of equality.
In one important respect, however, they behave differently:
as explained in \cite[Sections~4.3 and 5.1]{Altenkirch}, one can \emph{not} show that any two inhabitants $e,f$ of an equality type $a =_X b$ are equal---with 
their equality now being given by the iterated equality type $e =_{a=_Xb} f$. 

The lack of uniqueness of those terms has given rise to a new way of thinking about them and interpreting them into the world of mathematical objects.
Instead of interpreting them as (set-theoretic) equalities between $a$ and $b$ in the set interpreting $X$, 
one can interpret them as \emph{paths} from $a$ to $b$ in a space interpreting $X$.

This intuition is made formal in Voevodsky's \emph{simplicial set model} \cite{simplicial} which satisfies an additional interesting property:
given two types $X$ and $Y$, the interpretation of their equality type $X = Y$ is equivalent to the interpretation of their type of equivalences $X \simeq Y$ (see \cite[Section~5.4]{Altenkirch}).
This observation motivated Voevodsky to add this property as an axiom to Martin-L\"of type theory which he called the \emph{Univalence Axiom}.
The addition of the Univalence Axiom
turns Martin-L\"of type theory 
into univalent foundations. 

Obtaining an equivalence principle was one of the main motivations for Voevodsky in designing his univalent foundations:

\begin{quotation}
    [\ldots]
    My homotopy lambda calculus is an attempt to create a system which is  
    very good at dealing with equivalences. In particular it is supposed  
    to have the property that given any type expression $F(T)$ depending on  
    a term subexpression $t$ of type $T$ and an equivalence $t \to t'$ (a term of  
    the type $Eq(T;t,t')$) there is a mechanical way to create a new  
    expression $F'$ now depending on $t'$ and an equivalence between $F(T)$ and  
    $F'(T')$ (note that to get $F'$ one can not just substitute $t'$ for $t$ in $F $ 
    -- the resulting expression will most likely be syntactically  
    incorrect). [Email to Daniel R.~Grayson, Sept 2006]
 \end{quotation}

In the following sections, we describe how Voevodsky's goal is realized in univalent foundations.

\subsection{Indiscernibility of identicals in type theory}

In Martin-L\"of type theory (perhaps without the univalence axiom), identicals---that is, elements $x,y : T$ with an equality $e : x =_T y$ between them---are easily seen to be indiscernible. 
That is, for every type $T$ and $x,y : T$, we can find a function
\begin{equation}
              x =_T y \enspace \to \enspace  \forall \text{ properties } P, \left(P(x) \leftrightarrow P(y)\right)  \enspace .
\end{equation}
To better formulate this in the language of dependent type theory, (i) we will define this function for all $x, y$ at once, (ii) we will understand `properties  $P$' to be functions $P : T \to \Type$, and (iii) we replace the logical equivalence $P(x) \leftrightarrow P(y)$ with the type-theoretic equivalence
\[ \left( P(x) \simeq P(y) \right) := \left( \sum_{f:P(x) \to P(y) } \isEquiv(f) \right) \]
(where $\isEquiv$ is defined in Section~5.4 of Altenkirch's introduction \cite{Altenkirch}).

Our goal is hence to define a function
\begin{equation}            \transport: \prod_{(x,y:T)} \left( x =_T y  \enspace \to \enspace  \prod_{P: T \to \Type} \left(P(x) \simeq P(y)\right) \right) \enspace .
                   \label{eq:transport}
\end{equation}
To this end, recall that in order to define a map out of an equality type,
it suffices to define its image on $\refl_x : (x =_T x)$ for each $x: T$. Therefore, it suffices to show that there is a term 
\[            \transport(\refl): \prod_{(x:T)} \left(  \prod_{P: T \to \Type} \left(P(x) \simeq P(x)\right) \right) .\]
But then, for each $x : T$ and $P: T \to \Type$, we can set this to be the equivalence $1_{P(x)}: P(x) \simeq P(x)$ whose underlying function $P(x) \to P(x)$ is the identity function:
\[\transport(\refl)(x)(P) :\equiv 1_{P(x)}. \]

The function $\transport$ shows that any `property', or dependent type, $P: T \to \Type$ is invariant under equalities in $T$. In particular, given an equality $x =_T y$, we obtain functions $P(x) \to P(y)$ and $P(y) \to P(x)$ which allow us to transport terms of $P(x)$ or $P(y)$ back and forth along this equality.

\subsection{From equality to equivalence} \label{sec:etoe}

We have just seen that in Martin-L\"of type theory, identicals are indiscernible. Now we investigate how to expand this to get a full-blown equivalence principle from this fact. In short, we will see that in many circumstances, the equality type $x =_T y$ is itself equivalent to some structured equivalence appropriate for the type $T$. Then composing this equivalence with the $\transport$ function, we will obtain the equivalence principle \eqref{eq:equivalence_principle}.

To be precise, fix a type $T$. Given any notion $\sim_T$ of equivalence (or at least a reflexive relation) in a type $T$, we immediately obtain a function 
\begin{equation}   
	\idtoequiv : \prod_{x,y:T} \left(    (x = y) \xrightarrow{} (x \sim_T y) \right)
          \label{eq:eq_is_equiv}
\end{equation}
by setting $\idtoequiv(x,x)(\refl)$ to the reflexive term on $x$ in $x \sim x$ (since to define a function out of an equality type, it is enough to define it just at every occurence of $\refl$). 

Now we hope that for notions of equivalence $\sim_T$ already of interest to us, this function is actually an equivalence for all terms $x,y : T$, or more precisely, that the following type is inhabited.
\[ \prod_{x,y:T} \  \isEquiv \left(\idtoequiv_{x,y}\right)   \]
If this type is indeed inhabited, then for each $x,y : T$ we can take $\pi_1 (\isEquiv  (\idtoequiv(x,y)))$, the backwards function $(x \sim_T y) \to (x = y)$, and compose it with $\transport$ to obtain a function 
\[   \prod_{(x,y:T)} \left( x \sim_T y  \enspace \to \enspace  \prod_{P: T \to \Type} \left(P(x) \simeq P(y)\right) \right) \]
which is our equivalence principle.

Thus, in the next sections, we just aim to show that for certain types $T$ and notions of equivalence $\sim_T$, the function $\idtoequiv$ is indeed an equivalence.

\subsection{The univalence principle}

``Equality is equivalence for types'' is the slogan made precise by Voevodsky's 
univalence principle.
More precisely, the univalence principle asserts part of an equivalence principle for types:
it states that the canonical map 
\begin{align}\idtoequiv: \prod_{A,B:\Type} \left(    (A = B) \xrightarrow{} (A \simeq B) \right)
\label{eq:ua}\end{align}
 from equalities of types to equivalences of types is itself, for any types $A$ and $B$, an equivalence.
 Then, composing $\idtoequiv$ with $\mathsf{transport}$ as in the last section, we obtain an equivalence principle for types. 
 \[   \prod_{(A,B:\Type)} \left( A \simeq B  \enspace \to \enspace  \prod_{P: \Type \to \Type} \left(P(A) \simeq P(B)\right) \right) \]

The univalence principle is not provable in pure Martin-L\"of type theory  \cite{Martin-Lof-1972},
but needs to be postulated as an axiom---hence it is sometimes also called the ``univalence axiom''.
In extensions of Martin-L\"of type theory, as in the recently developed cubical type theory \cite{cohen_et_al:LIPIcs:2018:8475}, 
the univalence principle can be derived.

Building upon the equivalence principle for types---whether it is given as an axiom or as a theorem---one can
derive equivalence principles for other kinds of structures.
Establishing that $\idtoequiv$ is an equivalence
for other types and notions of equivalence is the subject of the next sections.

\section{The equivalence principle for set-level structures}

Now we turn our attention away from the type of all types and towards types of more specific mathematical objects. It turns out that for types of simple objects like propositions, sets, and monoids, the univalence axiom is enough to show the equivalence principle for these types' usual notion of equivalence. More precisely, in the presence of the univalence axiom, the function $\idtoequiv$ discussed in the last section is itself an equivalence. For an exploration and formalization of these ideas, see \cite{COQUAND20131105}.

\subsection{Propositions}

We call \emph{propositions} those types that have at most one inhabitant. We think of propositions as either being \emph{true} (when they are inhabited) or \emph{false} (when they are not inhabited). What should an equivalence of two propositions $P,Q$ be? Experience might indicate that such an equivalence should just be two functions
\[ f: P \leftrightarrows Q : g \]
so that $P$ is inhabited if and only if $Q$ is.
In fact, this notion of equivalence is the right one in the sense that it will validate the equivalence principle.

To be precise, we define
\begin{align*}
         &\mathsf{isProp} : \Type \to \Type
  \\
         &\mathsf{isProp}~A :\equiv \prod_{(x, y : A)} x = y
\end{align*}
whose inhabitants can be thought of as proofs that a type $A$ is a proposition.
A proposition is then a pair $(A,p)$ of a type $A$ and a proof $p : \mathsf{isProp}~A$, that is,
\[ \Prop :\equiv \sum_{A : \Type} \ \mathsf{isProp}~A \enspace .\]

Now we have, for $P \equiv (A,p)$ and $Q \equiv (B,q)$,
\begin{align}
            P = Q  \enspace & \simeq  \enspace  (A,p) = (B,q) 
                   \notag
          \\
                   \enspace & \simeq  \enspace  \sum_{(e : A = B)} (\transfib{\mathsf{isProp}}{e}{p}) = q
                   \label{eq:prop_sigma}
          \\	
                   \enspace & \simeq  \enspace  \sum_{(e : A = B)} 1
                   \label{eq:prop_isprop}
          \\	
                   \enspace & \simeq  \enspace  A = B
                   \notag
          \\	
                   \enspace & \simeq  \enspace   A \simeq B
                   \label{eq:prop_univalence}
          \\	
                   \enspace & \simeq  \enspace   (A \to  B) \times (B \to A)
                   \label{eq:prop_impl}
\end{align}
Equivalence~\eqref{eq:prop_sigma} above uses the fact that an equality between pairs is the same as pairs of equalities, 
where the second equality is ``heterogeneous'', i.\,e., requires a transport along the first equality to make it well-typed.
Equivalence~\eqref{eq:prop_isprop} uses the fact that being a proposition is itself a proposition, so that equality types between proofs of a proposition
are equivalent to the unit type.
Equivalence~\eqref{eq:prop_univalence} is given by the univalence principle, and equivalence~\eqref{eq:prop_impl} uses that $A$ and $B$ are propositions;
a pair of maps back and forth between types that are propositions automatically forms an equivalence of types.

Altogether, this means that $P$ and $Q$ are equal exactly if their underlying types $A$ and $B$ are logically equivalent---the expected notion of
equivalence for propositions.

\subsection{Sets}
We call \emph{sets} those types whose equality types are propositions,
\[ \Set :\equiv \sum_{A:\Type}\prod_{x,y:A}\mathsf{isProp}(x=y).  \]
Then a set in the type theory is a collection of terms, the equality types among which are either empty or contractible.

Given two sets $X,Y:\Set$, where $X=(A,p)$ and $Y=(B,q)$, the equivalences of types
\begin{align}  
  (X = Y) \enspace &\simeq \enspace (A \simeq B) 
   \label{eq:set_eq_equiv_eqv}
   \\& \simeq \enspace  \sum_{f : A \cong B} \mathsf{isCoherent}(f)
  \\ 
          \enspace &\simeq \enspace (A \cong B) 
   \label{eq:set_eqv_equiv_iso}
\end{align}
can be constructed.
Here, $A \cong B :\equiv \sum_{(f:A\to B)}\mathsf{isIso}f$ is the type of \emph{isomorphisms of types} between $A$ and $B$, 
and $\mathsf{isCoherent}(f)$ states an equality of equalities in $A$.
When $A$ is a set, the type $\mathsf{isCoherent}(f)$ is contractible (see the discussion in \cite[Section~5.4]{Altenkirch}),
and hence we obtain equivalence \eqref{eq:set_eqv_equiv_iso}.

\subsection{Monoids}
\label{sec:monoid}

The equivalence principle can be shown for many algebraic structures commonly encountered in mathematics, such as groups and rings.
Before presenting a general 
result to that extent in Section~\ref{subsec:categories}, in this section, we study in detail the case of monoids (which was formalized in \cite{COQUAND20131105}).
This particular case exemplifies many of the concepts and results used in general.

A \emph{monoid} is a tuple $(M,\mu,e,\alpha,\lambda,\rho)$ where
\begin{enumerate}
 \item $M : \Set$
 \item $\mu : M \times M \to M$ (multiplication)
 \item $e : 1 \to M$ (neutral element)
 \item $\alpha : \prod_{(a,b,c:M)} \mu(\mu(a,b),c) = \mu(a,\mu(b,c))$ (associativity)
 \item $\lambda : \prod_{(a:M)} \mu(e,a) = a$ (left neutrality)
 \item $\rho : \prod_{(a:M)} \mu(a,e) = a$ (right neutrality)
\end{enumerate}
Given two monoids $\mathbf{M}\equiv(M,\mu,e,\alpha,\lambda,\rho)$ and $\mathbf{M'}\equiv(M',\mu',e',\alpha',\lambda',\rho')$,
a \emph{monoid isomorphism} is a bijection $f : M \cong M'$ between the underlying sets that preserves multiplication and neutral element.
We can derive an equivalence between the equality type and the isomorphism type between any two monoids as follows:

  \begin{align}
      \id{\mathbf{M}}{\mathbf{M'}} 
        \enspace & \simeq  \enspace \id{(M,\mu,e)}{(M',\mu',e')} 
         \label{eq:monoid_prop}
        \\
                 & \simeq \enspace \sum_{p : \id{M}{M'}} (\id{\transfib{Y \mapsto (Y \times Y \to Y)}{p}{\mu}}{\mu'})  
                   \notag   
                      \\
                                                                                & \hspace{1.5cm}  \times(\id{\transfib{Y \mapsto (1 \to Y)}{p}{e}}{e'}) 
                                                                                \notag
        \\
                           & \simeq \enspace \sum_{f : {M}\cong{M'}} \bigl(\id{f \circ m \circ (f^{-1} \times f^{-1})}{m'}\bigr)  
                             \label{eq:monoid_univalence}
                           \\
  									       & \hspace{1.5cm}  \times(\id{f \circ e}{e'}) 
  									       \notag
  	\\
                           & \simeq \enspace    {\mathbf{M}}\cong{\mathbf{M'}} 
                           \notag
 \end{align}
Here, the equivalence of types \eqref{eq:monoid_prop} uses the fact that the axioms of a monoid (the types of $\alpha$, $\lambda$, and $\rho$)
are propositions (compare also to \eqref{eq:prop_isprop} above).
The equivalence \eqref{eq:monoid_univalence} uses the univalence principle for types in the first component,
replacing an equality of sets by a bijection.
This translates to replacing ``transport along the equality'' by ``conjugating by the bijection'' in the second component.

The equivalence of types constructed above, from left to right, is pointwise equal to the canonical map
\begin{equation}   
      \prod_{\mathbf{M}, \mathbf{M'}} \id{\mathbf{M}}{\mathbf{M'}}  \xrightarrow{}   {\mathbf{M}}\cong{\mathbf{M'}} \enspace   
      \label{eq:monoid_ep}
\end{equation}
defined by equality elimination,
which shows that the latter is an equivalence of types.
In other words, we have just proved 
the equivalence principle \eqref{eq:eq_is_equiv} for the equivalence ${\mathbf{M}}\cong{\mathbf{M'}}$.

\subsection{Univalent categories}\label{subsec:categories}

We have seen in the preceding sections that the types of propositions, sets, and monoids all have a certain nice property -- they validate the equivalence principle. However, it is natural to consider such objects as each belonging to a category. In this section, we discuss those categories whose objects validate the equivalence principle.

In Section \ref{sec:history}, we saw that in order to avoid mentioning equality of objects, we can define a
 \emph{category} $A$ to consist of
  \begin{enumerate}
  \item a type $A_0$ of \emph{objects};
  \item for each $a,b:A_0$, a set $A(a,b)$ of \emph{arrows} or \emph{morphisms};
  \item for each $a:A_0$, a morphism $1_a:A(a,a)$;
  \item for each $a,b,c:A_0$, a function of type
    \[  A(a,b) \times A(b,c) \to A(a,c) \]
    denoted by $(f,g) \mapsto f \cdot g$;
  \item for each $a,b:A_0$ and $f:A(a,b)$, we have $\ell_f : \id f {1_a\cdot f}$ and $r_f: \id f {f\cdot 1_b}$;
  \item for each $a,b,c,d:A_0$ and $f:A(a,b)$, $g:A(b,c)$, $h:A(c,d)$, we have $\alpha_{f,g,h}: \id {f\cdot (g\cdot h)}{(f\cdot g)\cdot h}$.
  \end{enumerate}
The reason for asking the types of arrows to be \emph{sets} rather than arbitrary types is so that these categories behave as classical categories (and not any kind of higher category) and, in particular, so that
 the axioms---which
state equalities between arrows---are propositions, meaning that we do not need to state higher coherence axioms.
There is prima facie no condition of that kind on the type of objects of the category. 
However, it will turn out that the objects of a \emph{univalent category} form a groupoid (meaning that all of its equality types form sets).

A morphism $f : A(a,b)$ of the category $A$ is an \emph{isomorphism} if there is a morphism $g : A(b,a)$ that is left and right inverse to $f$,
that is
\[ \mathsf{isIso}f :\equiv \sum_{g:A(b,a)} (f\cdot g = 1_a) \times (g \cdot f = 1_b) \enspace . \]
We call $\mathsf{Iso}(a,b) :\equiv \sum_{f:A(a,b)}\mathsf{isIso}f$ the type of isomorphims from $a$ to $b$,
and for any $a : A_0$ we have  $1_a : \mathsf{Iso}(a,a)$.
We can define a function
\begin{align} 
                   &\idtoiso : \prod_{x,y:A_0} (x = y) \to \mathsf{Iso}(x,y) 
                   \label{eq:uac}
\end{align}
by setting $\idtoiso_{x,x}(\refl_x)$ to $1_x$ for every $x: A_0$ just as we did to define 
$\idtoequiv$ in Section \ref{sec:etoe}.

Now we call the category $A$ \emph{univalent} if $\idtoiso_{x,y}$ is 
an equivalence of types for every $x,y:A_0$.
To see why the adjective \emph{univalent} is used, compare the function in Display~\eqref{eq:uac} above to the one in Display~\eqref{eq:ua} underlying the univalence principle. The univalence principle asserts that equality and equivalence of types are the same; here, we assert that equality and isomorphism of objects of a category are the same.

In asserting that a category $A$ is univalent, we assert that the equality types $a = b$ among its objects are equivalent to the sets $\mathsf{Iso}(a,b)$ of isomorphisms among its objects. Since the property of ``being a set'' itself obeys the equivalence principle for types the equality types $a = b$ are themselves sets.
When a type's equality types are sets, we call the type a \emph{groupoid}. 

A categorical equivalence between univalent categories $A$ and $B$ gives rise to an isomorphism between them---indeed, the type $A \simeq B$
of adjoint equivalences is equivalent to the type $A \cong B$ of isomorphisms of categories.

With a set-theoretic reading of the univalence condition in mind, one could think that only skeletal categories are univalent.
However, one should keep in mind that in type theory, the equality type $x = y$ between two objects of a category can---and 
often does---have more than one element.
Consequently, in type theory, a category being univalent usually signifies that its type of objects has many equalities.
This difference is witnessed by the many examples of univalent categories given below, most of which are not skeletal.

With these definitions in place, the composite equivalence of types shown in Displays~\eqref{eq:set_eq_equiv_eqv} - \eqref{eq:set_eqv_equiv_iso}
can be restated as ``the category of sets is univalent''.
Similarly, the result of Section~\ref{sec:monoid} can be restated
as ``the category of monoids is univalent''.

Many categories that arise naturally are univalent, in particular,
\begin{itemize}
 \item the category of sets;
 \item the categories of groups, rings, etc.;
 \item the functor category $[A,B]$ if the target category $B$ is;
 \item a preorder, seen as a category, exactly if it is anti-symmetric.
\end{itemize}

To extend our list of univalent categories to other algebraic structures beyond monoids, 
we could simply redo constructions similar to those for monoids,
for groups, rings, and other structures of interest.
However, in doing so, we would observe that we are doing the same reasoning over and over again.
For instance, looking back at monoids, we used that the category of sets is univalent to show that
the category of monoids is univalent, in step~\eqref{eq:monoid_univalence}.
This is due to the fact that ``monoids are sets with additional structure'',
and monoid isomorphisms are isomorphisms of sets preserving this structure.
Similarly, ``groups are monoids with additional structure'', 
and we would expect to reuse the equivalence of Display~\eqref{eq:monoid_ep} 
when building an equivalence between the equality types of groups on the one hand,
and of group isomorphisms on the other hand.
\emph{Displayed categories} as presented in \cite{ahrens_et_al:LIPIcs:2017:7722} are a convenient tool
for such modular reasoning about categories built step-by-step from simpler ones.
In particular, Proposition~43 and Theorem~44 of \cite{ahrens_et_al:LIPIcs:2017:7722} allow one
to show that a category built from a simpler one using the framework of displayed categories 
is univalent, provided the simpler one is univalent and the ``extra data'' making the 
difference between the two categories satisfies some condition.
That result validates the \emph{Structure Identity Principle} \cite[Theorem~9.8.2]{hottbook}.

\section{The equivalence principle for (higher) categorical structures}

We saw in the previous sections that for types of simple structures like propositions, sets, and monoids, the equivalence principle comes along with the univalence axiom. Now we see that for more complication structures, like categories, the equivalence principle only holds for certain well-behaved categories.

The most common notion of equivalence between two categories $A$ and $B$ is unsurprisingly called an equivalence $A \simeq B$. It consists of two functors $F: A \leftrightarrows B: G$ and natural isomorphisms $1_B \cong FG$ and $1_A \cong GF$ (see \cite{maclane}).
An equivalence of categories ``transports'' categorical structures, such as limits, between categories,
and is hence considered the right notion of sameness for categories in most contexts.
Can we show that the equality type between two categories, $A = B$, 
is the same as the type of categorical equivalences $A \simeq B$?
The answer is that while this is not the case for arbitrary categories,
it is the case when $A$ and $B$ are univalent.

For any two categories $A$ and $B$, the univalence axiom implies that the function from equalities
to isomorphisms (a stricter notion of sameness of categories) given by equality elimination is an equivalence \cite[Lemma~6.16]{rezk_completion}:
\begin{equation}   
	 (A = B) \xrightarrow{~~\simeq~~} (A \cong B) 
          \label{eq:eq_is_catiso}
\end{equation}
Furthermore, if $A$ and $B$ are univalent categories, then the type of isomorphisms between them is
equivalent to that of categorical equivalences \cite[Lemma~6.15]{rezk_completion}:
\begin{equation}   
	 (A \cong B) \xrightarrow{~~\simeq~~} (A \simeq B) 
          \label{eq:catiso_is_catequiv}
\end{equation}
Composing these two equivalences yields the desired equivalence of types between equalities and categorical equivalences.

The example of categories shows that, in order to obtain the equivalence principle for mathematical structures that naturally form 
bicategory, one needs to impose a ``univalence'' condition on those structures.
Defining such a univalence condition for general structures is the subject of active research.

\paragraph*{Acknowledgments}

We are very grateful to Deniz Sarikaya and Deborah Kant for their editorial work and their encouragement,
and to an anonymous referee for providing valuable feedback.
Furthermore, we would like to thank all the organizers of the FOMUS workshop---%
Balthasar Grabmayr, Deborah Kant, Lukas K\"uhne, Deniz Sarikaya, and Mira Viehst\"adt---%
for giving us the opportunity to discuss and compare different foundations of mathematics.

This material is based upon work supported by the Air Force Office of
Scientific Research under award numbers FA9550-16-1-0212 and FA9550-17-1-0363.

\bibliographystyle{plain}
\bibliography{literature}

\begin{thebibliography}{10}

\bibitem{rezk_completion}
Benedikt Ahrens, Krzysztof Kapulkin, and Michael Shulman.
\newblock {Univalent categories and the Rezk completion}.
\newblock {\em Mathematical Structures in Computer Science}, 25:1010--1039,
  2015.

\bibitem{ahrens_et_al:LIPIcs:2017:7722}
Benedikt Ahrens and Peter~LeFanu Lumsdaine.
\newblock Displayed categories \emph{(conference version)}.
\newblock In Dale Miller, editor, {\em 2nd International Conference on Formal
  Structures for Computation and Deduction (FSCD 2017)}, volume~84 of {\em
  Leibniz International Proceedings in Informatics (LIPIcs)}, pages 5:1--5:16.
  Leibniz-Zentrum f{\"u}r Informatik, 2017.

\bibitem{Altenkirch}
Thorsten Altenkirch.
\newblock Na{\"i}ve type theory.
\newblock In Stefania Centrone, Deborah Kant, and Deniz Sarikaya, editors, {\em
  Reflections on the Foundations of Mathematics: Univalent Foundations, Set
  Theory and General Thoughts}, pages 101--136. Springer International
  Publishing, Cham, 2019.

\bibitem{10.1093/philmat/nkt030}
Steve Awodey.
\newblock {Structuralism, Invariance, and Univalence}.
\newblock {\em Philosophia Mathematica}, 22(1):1--11, 10 2013.

\bibitem{MR539867}
Georges Blanc.
\newblock {\'Equivalence naturelle et formules logiques en th\'eorie des
  cat\'egories}.
\newblock {\em Arch. Math. Logik Grundlag.}, 19(3-4):131--137, 1978/79.

\bibitem{oberwolfach}
Samuel Buss, Ulrich Kohlenbach, and Michael Rathjen.
\newblock {Oberwolfach Reports -- Mathematical Logic: Proof Theory,
  Constructive Mathematics}.
\newblock pages 2963--3002.
\newblock \url{https://doi.org/10.4171/OWR/2011/52}.

\bibitem{cohen_et_al:LIPIcs:2018:8475}
Cyril Cohen, Thierry Coquand, Simon Huber, and Anders M{\"o}rtberg.
\newblock {Cubical Type Theory: A Constructive Interpretation of the Univalence
  Axiom}.
\newblock In Tarmo Uustalu, editor, {\em 21st International Conference on Types
  for Proofs and Programs (TYPES 2015)}, volume~69 of {\em Leibniz
  International Proceedings in Informatics (LIPIcs)}, pages 5:1--5:34,
  Dagstuhl, Germany, 2018. Schloss Dagstuhl--Leibniz-Zentrum fuer Informatik.

\bibitem{COQUAND20131105}
Thierry Coquand and Nils~Anders Danielsson.
\newblock Isomorphism is equality.
\newblock {\em Indagationes Mathematicae}, 24(4):1105 -- 1120, 2013.
\newblock In memory of N.G. (Dick) de Bruijn (1918–2012).

\bibitem{MR0412249}
Peter Freyd.
\newblock Properties invariant within equivalence types of categories.
\newblock In {\em Algebra, topology, and category theory (a collection of
  papers in honor of {S}amuel {E}ilenberg)}, pages 55--61. Academic Press, New
  York, 1976.

\bibitem{simplicial}
Chris Kapulkin and Peter~LeFanu Lumsdaine.
\newblock {The Simplicial Model of Univalent Foundations (after Voevodsky)}.
\newblock {\em J. Eur. Math. Soc.}
\newblock \href{https://arxiv.org/abs/1211.2851}{arXiv:1211.2851}.

\bibitem{leibniz}
Gottfried~Wilhelm Leibniz.
\newblock {\em Philosophical Papers and Letters}, volume~2 of {\em Synthese
  Historical Library}.
\newblock Springer-Verlag, 1989.

\bibitem{maclane}
Saunders Mac~Lane.
\newblock {\em Categories for the working mathematician}, volume~5 of {\em
  Graduate Texts in Mathematics}.
\newblock Springer-Verlag, New York, second edition, 1998.

\bibitem{MR1678360}
Michael Makkai.
\newblock Towards a categorical foundation of mathematics.
\newblock In {\em Logic {C}olloquium '95 ({H}aifa)}, volume~11 of {\em Lecture
  Notes Logic}, pages 153--190. Springer, Berlin, 1998.

\bibitem{Martin-Lof-1972}
Per Martin-L{\"o}f.
\newblock An intuitionistic theory of types.
\newblock In Giovanni Sambin and Jan~M. Smith, editors, {\em Twenty-five years
  of constructive type theory ({V}enice, 1995)}, volume~36 of {\em Oxford Logic
  Guides}, pages 127--172. Oxford University Press, 1998.

\bibitem{hottbook}
The {Univalent Foundations Program}.
\newblock {\em Homotopy Type Theory: Univalent Foundations of Mathematics}.
\newblock \url{http://homotopytypetheory.org/book}, Institute for Advanced
  Study, 2013.

\end{thebibliography}

\end{document}